\documentstyle[12pt]{article}
\input{amssym.def}
\input{amssym.tex}

\newcommand{\rest}{\restriction}
\newcommand{\ov}{\overline}
\newcommand{\oveta}{\ov \eta}
\newcommand{\overeta}{\ov \eta}
\newcommand{\ovmu}{\ov \mu}
\newcommand{\ovrho}{\ov \rho}
\newcommand{\overro}{\ov \rho}
\newcommand{\oversigma}{\ov \sigma}
\newcommand{\ovS}{\ov S}
\newcommand{\bbackslash}{\setminus}
\newcommand{\CE}{{\cal E}}
\newcommand{\fP}{{\cal P}}
\newcommand{\fc}{\frak{c}}
\newcommand{\dom}{{\rm dom}}
\newcommand{\range}{\mbox{\rm range}}
\newcommand{\code}{\mbox{\rm code}}
\newcommand{\Aone}{{\it A1}}
\newcommand{\At}{{\it A2}}
\newcommand{\Athree}{{\it A3}}
\newcommand{\Afour}{{\it A4}}
\newcommand{\Afive}{{\it A5}}

\newcommand{\qed}{\kern 5pt\vrule height8pt width6.5pt depth2pt}
\newcommand{\forces}{{\: \Vdash\:}}

\def\newtheorems{\newtheorem{theorem}{Theorem}[section]

                 \newtheorem{lemma}[theorem]{Lemma}
                 \newtheorem{claim}[theorem]{Claim}
                 
                 \newtheorem{definition}[theorem]{Definition}

                 }

\newtheorems

\title{Coding with ladders a well ordering of the reals}
\author
{Uri Abraham\\
Department of Mathematics and Computer Science, \\
Ben-Gurion University, Be'er-Sheva, Israel \\
 and\\
Saharon Shelah 
\thanks{This research was supported by The Israel Science Foundation
founded by the Israel Academy of Sciences and Humanities. Publication \# 485.}
\\
Institute of Mathematics\\
The Hebrew University, Jerusalem, Israel}
\begin{document}
\bibliographystyle{plain}
\maketitle
\begin{abstract}  Any model of ZFC + GCH has a generic extension
(made with a poset of size $\aleph_2$) in which the following hold: 
$MA + 2^{\aleph_0} = \aleph_2 +$ {\em there exists a}
 $\Delta^2_1$-{\em well ordering of
the reals.}  The proof consists in iterating posets designed to change at will
the guessing properties of ladder systems on $\omega_1$.  Therefore, the study
of such ladders is a main concern of this article. 
\end{abstract}

\section{Preface}
The character of possible well-orderings of the reals is a main
theme in set theory, and the work on long projective well-orderings
 by L. Harrington \cite{Harrington} can be cited as an example.
  There, the relative consistency of ZFC +
MA +$ 2^{\aleph_0} > \aleph_1$ with the existence of a $\Delta^1_3$
well-ordering of the reals is shown.  A different type of
 question is to ask about the  
impact of large cardinals on definable well-orderings.
Work of Shelah and Woodin \cite{ShelahWoodin}, and Woodin
\cite{Woodin}
is relevant to this type of question.  Assuming 
in $V$ a cardinal which is both measurable and Woodin, Woodin
\cite{Woodin} proved that 
if CH holds,
 then there is no $\Sigma^2_1$ well-ordering of the reals. 
 This result raises two questions:
\begin{enumerate}
  \item If large cardinals and  CH are assumed in $V$, can the 
 $\Sigma^2_1$ result be strengthen to $\Sigma^2_2$? 
That is, is there a proof that large cardinals and CH imply no $\Sigma^2_2$
well-orderings of the reals? 
\item  What happens if CH is not
assumed?
\end{enumerate}

Regarding the first question, Abraham and Shelah \cite{AS1} describes 
a poset of size 
$\aleph_2$ (assuming GCH) which generically adds no reals and 
provides a $\Delta^2_2$ well-ordering of the reals.
Thus,
if one starts with any universe with a large cardinal $\kappa$, one can
extend this universe with a small size forcing and obtain a $\Delta_2^2$
well-ordering of the reals.  Since small forcings  
will not alter the assumed largeness of a cardinal in $V$,
the answer to question 1 is negative.

Regarding the second question, Woodin (unpublished) uses an inaccessible
cardinal $\kappa$ to obtain a generic extension in which 
\begin{enumerate}
\item 
 MA for $\sigma$-centered posets +
$2^{\aleph_0} = \kappa$, and 
\item there is a $\Sigma^2_1$
well-ordering of the reals.  
\end{enumerate}
Solovay \cite{Solovay}
shows that the inaccessible cardinal is dispensable: any model of ZFC has a
small size forcing extension in which the following holds:
\begin{enumerate}
\item MA for $\sigma$-centered posets + 
$2^{\aleph_0} = \aleph_2$, and 
\item there is a $\Sigma^2_1$ well-ordering of
the reals. 
\end{enumerate}

In \cite{AS2} we show how Woodin's result can be strengthened to obtain the
full Martin's axiom. We prove there that if $V$ satisfies the GCH and
 contains an inaccessible cardinal $\kappa$, 
then there is a poset of cardinality $\kappa$ that
gives generic extensions in which 
\begin{enumerate}
\item MA + $2^{\aleph_0} = \kappa$, and 
\item there is a $\Sigma^2_1$ well-ordering of
the reals. 
\end{enumerate}

Our aim in this paper
 is to show that the inaccessible cardinal is not really necessary,
even to get the full Martin's Axiom.
\begin{theorem} 
\label{main}
Assume $2^{\aleph_0} = \aleph_1$ and $2^{\aleph_1} = \aleph_2$.
 There is a forcing poset of size $\aleph_2$ that provides a cardinal
preserving extension in which Martin's Axiom
 $+2^{\aleph_0} = \aleph_2$ holds, and there is a $\Sigma^2_1$ well-ordering of
the reals.  In fact, there is even a $\Sigma^{2[\aleph_1]} $ well-ordering of
the reals there.
\end{theorem}
The concepts $\Sigma^2_1$ and $\Sigma^{2[\aleph_1]} $ will soon be defined,
but first we shall point to what we consider to be the main novelty of this
paper, the use of ladder systems as coding devices.  A ladder over $S
\subseteq \omega_1$ is a sequence $\overline{\eta}=
\langle \eta_\delta \mid \delta \in S\rangle$
 where
$\eta_\delta : \omega \to \delta$ is increasing and 
 cofinal in $\delta$.  Two ladders over
$S$, $\overline{\eta}'$ a subladder of $\overline{\eta}$,
 may encode a real (a subset of $\omega$). 
Namely the coding of a real $r$ is expressed by the relationship
 between $\eta'_\delta$ and $\eta_\delta$  (for every $\delta$).
Splitting $\omega_1$ into $\aleph_2$ pairwise almost disjoint stationary sets,
it is possible to encode $\aleph_2$ many reals (and hence a well-ordering)
using $\aleph_2$ pairs of ladder sequences.  Of course, we need some property
that ensures uniqueness of these ladders, in order to make this well-ordering
definable.  Such a property will be obtained in relation with the guessing
power of the ladders.  A ladder system $\langle \eta_\delta \mid \delta
 \in S\rangle$ is said to be club (closed unbounded set)
 guessing if for every closed unbounded $C
\subseteq \omega_1,\;[\eta_\delta] \subseteq^* C$ for some $\delta \in S$.  It
turns out that there is much freedom to manipulate the guessing properties of
ladders, and, technically speaking, this shall be a main concern of the paper.

We now define the 
 $\Sigma^2_1$ and $\Sigma^{2[\aleph_1]} $ relations.
The structure with the membership relation on the 
collection of all hereditarily countable sets
 is denoted $H(\aleph_1)$.  Second-order formulas
over $H(\aleph_1)$ that contain $n$ alternations of quantifiers are denoted
$\Sigma^2_n$ when the external quantifier is an existential class  quantifier. 
Thus a $\Sigma^2_n$ formula has the form $$\exists X_1 \forall X_2\ldots X_n
\varphi (X_1,\ldots,X_n)$$ where $\varphi$ may only contain first-order
quantifiers over $H(\aleph_1)$ and predicates $X_1,\ldots,X_n$ are
interpreted as subsets of $H(\aleph_1)$.
(One can either write $X_i(s)$ treating $X_i$ as a predicate, or $s \in
X_i$ treating $X_i$ as a class.)
$\Sigma^2$ denotes the union of all $\Sigma^2_n$ formulas.

If the second-order quantifiers only quantify classes (subsets of 
$H(\aleph_1)$) of cardinality $\leq  \aleph_1$, 
then the resulting set of formulas is denoted 
$\Sigma^{2[\aleph_1]}_n$. 
So $\Sigma^{2[\aleph_1]}_1$ for example denotes second order formulas of
the form ``there exists a subset $X$ of $H(\aleph_1)$ of size $\leq
\aleph_1$ such that $\varphi(X)$'' where $\varphi$ is a first order formula.
 We write $\Sigma^{2[\aleph_1]}$, without a subscript,
for $\bigcup_{n<\omega} \Sigma^{2[\aleph_1]}_n$.

In Theorem \ref{main} above, 
we get a well-ordering which is $\Sigma^{2[\aleph_1]}$, and we
will explain now why $MA + 2^{\aleph_0} > \aleph_1$ implies that such a
relation is necessarily $\Sigma^2_1$.  This transformation
which replaces any number of quantifiers over sets of size $\aleph_1$ with a
 single existential quantifier over arbitrary subsets of $H(\aleph_1)$
 is a trick of Solovay's that was used by him in \cite{Solovay}.
  The basic idea is to use the almost-disjoint-sets
coding (Jensen and Solovay \cite{JensenSol}) 
in a way which will be sketched here.
\begin{theorem}[Solovay]
Assume MA+$2^{\aleph_0}> \aleph_1$.
  Any  $\Sigma^{2[\aleph_1]}$ formula $\varphi(\overline{x})$ 
 over $H(\aleph_1)$,
with free variables $x_1,\ldots,x_n$, is equivalent to a 
$\Sigma^2_1$ formula $\psi(\overline{x})$.
\end{theorem}
{\bf Proof}. It seems easier to prove first that every
 $\Sigma^{2[<{\fc}]}$
formula is equivalent with a $\Sigma^2_1$ formula. 
(The $\Sigma^{2[<{\fc}]}$
formulas are second order formulas over $H(\aleph_1)$ in which class quantification
occurs only for subset of $H(\aleph_1)$ of size less than continuum.)
Then the theorem follows because the
$\Sigma^{2[\aleph_1]}$ classes are a naturally characterized subclass of the
 $\Sigma^{2[<{\fc]}}$.

So let $\varphi(x)$ be any $\Sigma^{2[<\fc]}$ formula.
  The equivalent $\Sigma^2_1$
formula $\psi$ begins as follows (with existential class quantifiers mixed
with first-order quantifiers which do not change the complexity of the
formula):  
\begin{quotation}
{\em There is a set
$\tau \subset \fP (\omega)$ such that the relation
 $$x <_\tau y\  \mbox{\it iff } y \setminus x \ \mbox{\it  is finite}$$
is a well-order of $\tau$
 such that there is no infinite $a \subseteq \omega$ with 
$a \subseteq^* x$ for all $x \in \tau$.
There is also a map $\mu :\tau \longrightarrow H(\aleph_1)$,
 which is {\em onto} $H(\aleph_1)$, and there is 
 a map $\rho : \tau \longrightarrow [\omega]^{\aleph_0}$ such that for distinct 
$x,y \in
\tau$,  $\tau(x)$ and $\tau(y)$ are
almost disjoint. ($ [\omega]^{\aleph_0}$ is the collection of infinite
subsets of $\omega$.)}
\end{quotation}
  Then $\psi$ continues 
 with first-order quantifiers
that replace the 
$\Sigma^{2[< \fc]}$ quantifiers of $\varphi$
 in the following manner.  To represent any $X \subseteq
H(\aleph_1)$ of size $< \fc$, look at the set $\mu^{-1} X \subseteq \tau$. 
Since its size is $<\fc$, there is by Martin's Axiom an infinite set
$a\subseteq \omega$ almost included in every
 set in $\mu^{-1}X$. Hence $\mu^{-1}X$ is bounded in $\tau$. 
 So there is 
$t_0$ in $\tau$ so that $\mu^{-1}(x) <_\tau t_0$ for every $x \in X$.  Now look at
the collection $\{ \rho(t)  \mid t <_\tau
 t_0\}$ of almost-disjoint sets (its cardinality
is $< \fc$) and use Martin's Axiom to
 encode with one $r$ the set $\rho[\mu^{-1}X]$.
That is find $r \subset \omega$ such that for $t <_\tau t_0$, $\rho(t) \cap
r$ is finite iff $\mu(t) \in X$. Then $r$ and $t_0$  represent $X$.

\section{Ladder systems}

 The notation $A \subseteq^* B$ is used for ``almost inclusion'' on subsets
of $\omega_1$,
meaning that $A\bbackslash B$ is  finite.  Similarly $A =^* B$ is defined
if $A\subseteq^* B$ and $B \subseteq^* A$.
$A \not =^* B$ is the negation of $A=^* B$.

\begin{definition}
\label{Dladder}
\begin{enumerate}
\item 
A {\em ladder system} over $S \subseteq \omega_1$ (consisting of limit ordinals)
 is a sequence $\overline{\eta} = \langle \eta_\delta \mid \delta \in S\rangle$,
where $\eta_\delta$ is an increasing $\omega$-sequence converging to $\delta$.
 $S$ is called ``the domain'' of $\overline{\eta}$, and is denoted
$\dom(\oveta)$. $\oveta$ is called ``trivial'' if $\dom(\oveta)$ is
non-stationary.
The range of $\eta_\delta$ is denoted $[\eta_\delta]$ (so
$[\eta_\delta]=\{\eta_\delta(i) \mid i \in \omega \}$),
and $\bigcup_{\delta\in S}$
$[\eta_\delta]$ is the ``range'' of $\oveta$. 
 So, $[\eta_\delta] \subset^* C$
means that, except for finitely many $k$'s, $\eta_\delta(k) \in C$ always
holds.  
\item 
Let $\overline{\eta}$ and $\overline{\mu}$ be two ladder systems.  We
say that $\overline{\eta}$ and $\overline{\mu}$ are {\em almost disjoint}
 iff, for some club $C \subseteq \omega_1$,
 for any $\delta \in C \cap \dom (\overline{\eta}) \cap \dom
(\overline{\mu}),\; [\eta_\delta] \cap [\mu_\delta] =^* \emptyset$.  
\item
We say that  $\overline{\eta}$ is a subladder of $\overline{\mu}$ iff 
the following holds
for some club $C \subseteq \omega_1$: 
$$C \cap \dom(
\overline{\eta})\subseteq \dom( \overline{\mu})\mbox{ and for }\delta \in  
C \cap \dom (\overline{\eta}), \;[\eta_\delta] \subseteq^* [\mu_\delta].$$
In such a case we
 write $\overline{\eta}\lhd \overline{\mu}$.  Also, $\overline{\eta} =^* 
\overline{\mu}$ iff both $\overline{\eta} \lhd \overline{\mu}$ and 
$\overline{\mu} \lhd \overline{\eta}$.
That is,
$\overline{\eta} =^* \overline{\mu}$ iff there is a club
set $C \subseteq \omega_1$ such that $\dom(\overline{\eta}) \cap C =
 \dom(\overline{\mu}) \cap C$, and $[ \eta_\delta] =^* [ \mu_\delta]$
 for $\delta
\in \dom (\overline{\eta}) \cap C$.

\item
The difference ladder $\overline{\rho} = \overline{\eta}\bbackslash \overline{\mu}$ is
defined by
\[ [\rho_\delta] = \left\{ \begin{array}{ll} 
[\eta_\delta]& \mbox{if} \;\delta \in \;
\dom(\overline{\eta})\bbackslash\;
\dom (\overline{\mu}) \\
{[ \eta_\delta ]} \bbackslash[ \mu_\delta]& \;\mbox{if this set is infinite} \\
\mbox{undefined}& \mbox{otherwise}\end{array} \right. \]
It is the $\lhd$-maximal ladder included in $\overline{\eta}$ and
(almost) 
disjoint
from $\overline{\mu}$.
\item
Given any $A \subseteq \omega_1$, the restriction ladder  $\overline{\eta}
 \rest
 A$ is naturally defined, and its domain is $A \cap\dom(\overline{\eta})$.  If $x \subseteq \omega$ is infinite, then 
  $\overline{\eta} \rest x$ means something else:  it is obtained
 by enumerating $x = \{x_k\mid k \in \omega\}$ in increasing order, and setting 
$(\overline{\eta} \rest x) = \overline{\rho}$ where 
$\rho_\delta(k) = \eta_\delta(x_k)$ for every $\delta \in
\dom(\overline{\eta})$.
\end{enumerate}
\end{definition}

We shall define some properties of ladders (in fact, of $=^*$ equivalence
classes).
\begin{definition}
 Let $\overline{\eta}$ be a ladder over $S$.
\begin{enumerate}
\item We say that $\ov{\eta}$ is {\em club-guessing}
 iff for every club $C \subseteq \omega_1$ there is $\delta \in S$ such that
$[ \eta_\delta] \subseteq^* C$.  (So, in this case, $\delta \in C$, and
hence  $\dom(\overline{\eta})$ is stationary if $\overline{\eta}$ is club
guessing.)  
For brevity, we may use the term {\em guessing} instead of club-guessing.
\item We say that $\overline{\eta}$ is {\em strongly} 
club guessing (or just strongly guessing) 
iff for any club $C \subseteq \omega_1$ for some club $D$, if 
$\delta \in D \cap S$ then $[ \eta_\delta] \subseteq^* C$.  If $\overline{\eta}$
 is strongly guessing and $\overline{\rho} \lhd \overline{\eta}$, then clearly 
$\overline{\rho}$ is also strongly guessing.  (Be careful: if
 $\overline{\rho} \lhd \overline{\eta}$ and $\overline{\eta}$ is guessing,
 you cannot infer that $ \overline{\rho}$ is guessing, unless $
 \overline{\rho}$ is non-trivial.)   
The trivial ladder is (trivially) strongly guessing, and hence we cannot say
that a strongly guessing ladder is always guessing. A strongly guessing
non-trivial ladder is, of course, guessing.

\item  We say that a club set
 $C \subseteq \omega_1$ {\em avoids} $\overline{\eta}$ iff for every 
$\delta \in  S$ (except a non-stationary set), 
$[ \eta_\delta] \cap C =^* \emptyset$. 
\item We say that
 $\overline{\eta}$ is {\em avoidable} iff some club set avoids 
$\overline{\eta}$. If every ladder
over $S$ is avoidable, then we say that $S$ itself is {\em avoidable}.  Hence,
in particular, if $S$ is non-stationary, then $S$ is avoidable.  
Remark that if $\oveta$ is avoidable, then 
$\oveta$ is non-guessing.
So $\oveta$ is strongly guessing and avoidable iff $\oveta$ is trivial.
The collection of all avoidable sets forms an ideal which will be shown to be
normal in the following subsection. 
\item
{\bf Maximal ladders}.  Suppose that $\overline{\eta}$ is  some strongly 
guessing ladder over $S$, and $X \supseteq S$ is a subset of
$\omega_1$.  If every ladder
over $X$ and (almost) 
disjoint from $\overline{\eta}$ is avoidable, then we say that 
$\overline{\eta}$ is {\em maximal for} $X$.  In such a case, for every $X'
\subseteq X,\;\overline{\eta} \rest X'$ is 
maximal for $X'$.  The
trivial ladder $\emptyset$ is trivially maximal for any avoidable set.  Our
terminology may be misleading because a maximal ladder for $X$ is not necessarily
defined over $X$, it is rather the maximality which is for $X$.  Thus, if
 $\overline{\eta}$ is maximal for $X$, then $\overline{\mu} \lhd
 \overline{\eta}$ for every strongly guessing ladder $\overline{\mu}$ over a
subset of $X$. 
(Because $\ovmu \setminus \oveta$ is, in that case, strongly guessing and
disjoint from $\oveta$, and is hence avoidable. Thus $\dom(\ovmu \setminus
\oveta)$ is not stationary, and hence $\ovmu \lhd \oveta$.)
 Hence if both $\overline{\mu}$ and $\overline{\eta}$ are
maximal for $X$, then $\overline{\mu} =^* \overline{\eta}$.  We denote this
unique ladder, maximal for $X$, by $\chi(X)$.
\end{enumerate}
\end{definition}
It is easy to see that if $\overline{\eta}$ is maximal for $X$ and $X_0
\subseteq X$ then $\overline{\eta}\rest X_0$ is maximal for $X_0$.
\subsection{Ideals connected with ladders}

 We are going to define four ideals on $\omega_1$:
 the ideal of
non-guessing restrictions, denoted $I_{\overline{\eta}}$, 
  the ideal of avoidable sets, denoted $I_0$, the ideal of
maximal guesses, denoted $I_1$, and the ideal of bounded intersections,
$I(\overline{S})$. Then we will prove that all are normal ideals.

\begin{definition}
\label{Dideals}
\begin{description}
\item
[The ideal of non-guessing restrictions.]  Let $\overline{\eta}$ be a
guessing ladder over $X$.  The collection of all
 subsets $S \subseteq \omega_1$ for
which $\overline{\eta} \rest S$ is not guessing is a proper, normal
ideal, denoted $I_{\overline{\eta}}$. 
\item
 [The ideal of avoidable sets.]  $S \in I_0$ iff every ladder system over
$S$ is avoidable.
\item
 [The ideal of maximal guesses.]  The ideal $I_1$ is the collection of all
sets $X \subseteq \omega_1$ such that there is a maximal ladder for $X$.

So $S \in I_1$ iff there is a strongly guessing ladder system
$\overline{\eta}$ such that $\dom(\overline{\eta}) \subseteq S$ and any ladder
over $S$ and disjoint from $\overline{\eta}$ is avoidable.  As said above,
this unique ladder $\overline{\eta}$ is denoted $\chi(S)$.  (Uniqueness is up
to $=^*$, where non-stationary sets and finite differences do not count.)

In case $S \in I_0$, then $S \in I_1$, and $\chi(S)$ is the trivial (empty)
ladder $\emptyset$.
So \begin{equation} \label{e0}
I_0 \subseteq I_1.\end{equation}

\item
 [The ideal of bounded intersections.]  Let 
 $\overline{S} = \langle S_i \mid  i \in \omega_2\rangle$ be a collection of
$\aleph_2$ 
 stationary subsets of $\omega_1$ such that the intersection
of any two is non-stationary (we say that $\overline{S}$ is a sequence
of pairwise almost disjoint stationary sets).  The ideal 
$I(\overline{S})$ consists of those sets $H \subseteq \omega_1$ for which 
$$ |\{ i \in \omega_2 \mid H \cap S_i \mbox{ is stationary }\}| \leq
\aleph_1.$$  Sets in $I(\overline{S})$ will also be called 
$\overline{S}$-small sets. 
 It may seem that $I(\overline{S})$ is not connected to ladders, but
we will later show the consistency of $I(\overline{S}) = I_1$.  
\end{description}
\end{definition}
\begin{lemma}
All four ideals are normal.
\end{lemma}
{\bf Proof}.
An ideal on $\omega_1$ is said to be normal if it is closed under diagonal
unions. 
\begin{description}
\item[The ideal of non-guessing restrictions.]
Let $\overline{\eta}$ be a guessing ladder over $X$. 
 To prove normality of $I_{\oveta}$,
 suppose $A_\xi \in
 I_{\overline{\eta}}$, for $\xi < \omega_1$.  Thus, 
for every $\xi$ there is a club
set $C_\xi$ such that $\delta \in  A_\xi \cap X \Rightarrow [\eta_\delta]
\not\subset^* C_\xi$.  Let $$ A = \nabla_{\xi\in\omega_1} A_\xi
 \stackrel{\mbox{def}}{=}
\{ \alpha \in \omega_1 \mid  \exists \xi < \alpha (\alpha\in A_\xi)\}$$ 
be the diagonal
union, and $C = \Delta_{\xi\in\omega_1} C_\xi$ be the diagonal intersection
 of the club sets.  Then $A \in I_{\ov{\eta}}$ because for $\delta \in A
\cap X$, $ [\eta_\delta] \not\subset^* C$.

\item[The ideal of avoidable sets.]
We check that $I_0$ is normal.  Suppose $S_\xi \in I_0$ for $\xi \in
\omega_1$, and let $S = \nabla_{\xi\in\omega_1} S_\xi$ be the
diagonal union.
 Let 
$\overline{\eta} = \langle \eta_\delta \mid \delta \in S\rangle$ be any ladder over
$S$, and we will show that $\overline{\eta}$
is avoidable and hence that $S \in I_0$. 
Indeed, a slightly more general fact will be used later: 
\begin{quotation}
\noindent
If $S_\xi \subseteq
\omega_1$ are arbitrary sets, $S =\nabla_{\xi\in\omega_1}
S_\xi$, and $\overline{\eta}$ is a ladder over $S$ such that 
$\overline{\eta} \rest
 S_\xi$ is avoidable for every $\xi \in \omega_1$,
then $\overline{\eta}$ is avoidable.
\end{quotation}
 To see this, let $C_\xi$ for $\xi \in \omega_1$ be  a club set that avoids 
$\overline{\eta} \rest S_\xi$, and  let $C = \Delta_{\xi\in\omega_1}
C_\xi$ be their diagonal intersection.  Then $C$ avoids
 $\overline{\eta}$, as can easily be checked.
\item[The ideal of maximal guesses.]
We prove that $I_1$ is normal.  So suppose that  $S_\xi \in I_1$ for $\xi \in
\omega_1$ are given, and $S = \nabla_{\xi\in\omega_1} S_\xi$
is their diagonal
union.  We must prove that $S \in I_1$.  
First we claim that the sets $\{ S_\xi \mid \xi \in \omega_1 \}$ may be
assumed to be pairwise disjoint. Indeed, define $S^*_\xi =
S_\xi \setminus \bigcup \{  S_{\xi'} \mid \xi' < \xi \}$.
 Then $S =
\nabla S^*_\xi$, and the sets  $S^*_\xi$ 
are in $I_1$ and are pairwise
disjoint. So we do assume now that the $S_\xi$'s are pairwise disjoint.
For every $\xi \in
\omega_1,\;\chi(S_\xi)$ is a strongly guessing ladder over its domain 
$S^0_\xi \subseteq
S_\xi$ (and $S^0_\xi = \emptyset$ when $S_\xi \in I_0$).  Define 
$S^0 = \nabla_{\xi\in\omega_1} S^0_\xi$. Clearly $S^0 \subseteq
S$. For $\delta \in S^0$
define $\eta_\delta$ to be $\chi(S_\xi)_\delta$ for the (unique) $\xi < \delta$
such that $\delta \in S^0_\xi$.

\noindent
{\bf Claim}:  $\overline{\eta} = \langle \eta_\delta \mid \delta \in S^0\rangle$
is maximal for $S$, and hence $S \in I_1$.

\noindent
{\bf Proof}.  We first prove that $\overline{\eta}$ is strongly guessing.
 Well, if $C \subseteq \omega_1$ is club, find for each $\xi \in \omega_1$ a
club set $D_\xi$ such that for $\delta \in D_\xi \cap
S^0_\xi,\; (\chi(S_\xi))_\delta \subseteq^* C$.  Now define $D =
\Delta_{\xi\in\omega_1} D_\xi$ to be the diagonal intersection.  It
follows that for every $\delta \in S^0 \cap D,\;[\eta_\delta] \subset^* C$.

To prove maximality, assume $\overline{\mu}$ is defined on $S$ and is disjoint
from $\overline{\eta}$. Then $\overline{\mu} \rest S_\xi$ is disjoint from
$\chi(S_\xi) \rest S_\xi \setminus (\xi+1)$. Hence $\overline{\mu}\rest
S_\xi$ is avoidable for every $\xi < \omega_1$, and by the proof of
normality of $I_0$, $\overline{\mu}$ is avoidable.

\item[The ideal of bounded intersections.]
Let $\overline{S} = \langle S_i \mid i \in \omega_2 \rangle$ be a
collection of pairwise almost disjoint stationary subsets of $\omega_1$
defining $I(\overline{S})$. If $H_\xi$ for $\xi \in \omega_1$ are in
$I(\overline{S})$, then there is a bound $j_0 < \omega_2$ such that for
every $j_0 \leq j < \omega_2$ $S_j \cap H_\xi$ is non-stationary. Hence
$\nabla_\xi S_j \cap H_\xi= S_j \cap \nabla H_\xi$ is
non-stationary, and thus $\nabla_\xi H_\xi\in I(\overline{S})$.
\end{description} 

\subsection{$A(\overline{S},\overline{\eta})$}
In this subsection we 
formulate a statement, $A(\overline{S},\overline{\eta})$, and show
that it implies $I(\overline{S}) = I_1$.  The consistency of
 $A(\overline{S},\overline{\eta})$ will be proved in the subsequent sections.

\begin{definition}
\label{DA}
 $A(\overline{S},\overline{\eta})$ is the conjunction of the
following six statements:
\begin{description}
\item{\Aone} $\overline{S} = \langle S_i \mid  i \in \omega_2\rangle$ is a sequence
 of pairwise almost disjoint stationary subsets of $\omega_1$.  
$\overline{\eta}$ is a ladder system, and  $\bigcup_{i<\omega_2} S_i \subseteq
\dom (\oveta)$.

\item{\At} Every ladder disjoint from $\overline{\eta}$ is avoidable.  (It 
immediately follows that if $\overline{\mu}$ is strongly guessing, then 
$\overline{\mu}\bbackslash \overline{\eta}$ is both avoidable and strongly
guessing and thus
$\overline{\mu}\bbackslash \overline{\eta} =^* \emptyset$, so that
 $\overline{\mu} \lhd \overline{\eta}$.)

\item{\Athree} For every $i < \omega_2,\;S_i \in I_1$.  In fact, $\chi(S_i)$ is
defined over $S_i$ (and it is a non-trivial strongly guessing ladder over
$S_i$ such that any ladder over a subset of
$S_i$ and disjoint from $\chi(S_i)$ is avoidable).
It follows by {\At} that $\chi(S_i) \lhd \overline{\eta}$.

\item{\Afour} If $X \subseteq \omega_1$ is such that $X \cap S_i$ is
non-stationary for every $i < \omega_2$, then $X$ is avoidable (equivalently, 
in view of (\At), 
$\overline{\eta} \rest X$ is avoidable).

\item{\Afive} If $X \subseteq \omega_1$ is not $\overline{S}$-small,
 $\overline{\rho}$ is a ladder over $X$ and $\overline{\rho} \lhd 
\overline{\eta}$, then there exists $i < \omega_2$  such that $S_i \subseteq
X$ and $\chi(S_i) \lhd \overline{\rho}$.

\item{A6} For every $i \in \omega_2$ either $(\chi(S_i),
\overline{\eta}\rest S_i)$ is clearly not encoding, or else $r =
d(\chi(S_i), \overline{\eta}\rest S_i)$ is defined, and in this case
 $r = d(\chi(S_j), \overline{\eta}\rest S_j)$ 
for unboundedly many $j$'s.
The meaning of this statement is clarified later in this subsection.

\end{description}
\end{definition}
$A'(\overline{T},\overline{\eta})$ is the following statement:  $\overline{T}
 = \langle T_i \mid i < \omega_2 \rangle$ is a sequence of pairwise almost
disjoint stationary subsets of $\omega_1$.  For every $i < \omega_2,\;T_i \in
I_1$, and if $S_i$ denotes $\dom(\chi(T_i))$, then
 $A(\overline{S},\overline{\eta})$ holds for $\overline{S} = \langle
 S_i \mid i < 
\omega_2\rangle$.

We first collect some simple consequences of
the first five statements of $A(\overline{S},\overline{\eta})$.
\begin{lemma}
\label{l21}
The first five statements of $A(\ov{S},\oveta)$ imply that:
\begin{enumerate}
\item
 If $\overline{\rho} \lhd \overline{\eta}$
is avoidable, then $\dom(\overline{\rho})$ is $\overline{S}$-small.  
\item
$I_0 \subseteq I(\overline{S})$.
\item If $\ovmu \lhd \oveta$ is strongly guessing, then $\dom(\ovmu)$ is
$\ovS$-small.
\item Actually: If $\overline{\mu}$ is strongly guessing, then $\dom(
\overline{\mu})$ is  $\overline{S}$-small.
\item $I_1 = I( \overline{S})$.
\end{enumerate}
\end{lemma}
{\bf Proof}. To prove {\it 1}, assume $\ovrho \lhd \oveta$ but
$X=\dom (\ovrho)$ is not $\ovS$-small. Then (\Afive) implies that, for
some $i < \omega_2$, $\chi(S_i) \lhd \ovrho$. Hence $\ovrho$ is not
avoidable (by (\Athree) which says that $\chi(S_i)$ is (strongly) guessing).

We prove {\it 2}. If $X\in I_0$ ($X$ is avoidable) then any ladder system over
$X$, and in particular $\overline{\eta}\rest X$, is avoidable. Hence (by
item {\it 1}) $\dom(\overline{\eta}\rest X)$ is $\overline{S}$-small. Thus $X$
is  $\overline{S}$-small (because $X=X_0 \cup X_1$ where $X_0 = X \cap
\bigcup_i S_i$ and $X_1 = X \setminus X_0$. $X_1$ is clearly
$\overline{S}$-small, and $X_0 = \dom (\overline{\eta}\rest X)$).

To prove {\it 3}, assume that $\dom(\overline{\mu})$ is not
$\overline{S}$-small. 
Split $\overline{\mu}$ into $\overline{\mu}^1$ and
 $\overline{\mu}^2$, two ``halves'' defined by taking $(\mu^1)_\delta$ to be
an infinite co-infinite subset of $\mu_\delta$ (for every $\delta \in
 \dom(\overline{\mu}))$, and letting $\overline{\mu}^2  = 
\overline{\mu}\bbackslash \overline{\mu}^1$.  If $X = \dom(\overline{\mu})$
is {\em not} $\overline{S}$-small, then, by (\Afive) applied to $\overline{\mu}^1$,
there is $i$ such that $S_i \subseteq X$ and
\begin{equation}
\label{e1}
\chi(S_i) \lhd \overline{\mu}^1.
\end{equation}
 Since $\overline{\mu}$ is strongly guessing,
 $\overline{\mu}^2 \rest S_i$ is strongly guessing (and
non-trivial as its domain is the stationary set $S_i$), but formula
(\ref{e1})
shows that $\ov{\mu}^2 \rest S_i$ is disjoint from $\chi(S_i)$, and this
contradicts the maximality of $\chi(S_i)$ for $S_i$.

  To prove {\it 4}, suppose that $\overline{\rho}$ is a strongly guessing
 ladder over $X$.  To show that $X \in I(\overline{S})$, we reduce this
claim to the case that $\overline{\rho}\lhd \overline{\eta}$.  Look at 
$\overline{\rho}\bbackslash \overline{\eta}$ and its domain
\[ X_1 = \{\delta \in X  \mid
[\rho_\delta] \setminus [\eta_\delta]
\;\mbox{is infinite} \} . \]
By (\At), $\overline{\rho}\bbackslash \overline{\eta}$ is avoidable.  But, as 
$\overline{\rho}$ is strongly guessing, any subladder of $ \overline{\rho}$ is
also strongly guessing, and hence $\overline{\rho}\bbackslash \overline{\eta}$
is strongly guessing and avoidable, which could only be if $X_1$ is
non-stationary.

Now set $X_2 = X\bbackslash X_1$, and $\overline{\mu} = \overline{\rho}
\rest X_2$.  Then $\overline{\mu}\lhd \overline{\eta}$ is strongly
guessing, 
 and hence by the previous item 
  $\dom(\overline{\mu})$ is $\overline{S}$-small.

Finally we prove {\it 5}.
 If $X \in I_1$ then $X =X_0 \cup X_1$, where $X_0
\in I_0$ and $X_1$ is the domain of a strongly guessing ladder---namely
$\chi(X)$. Hence $X\in I(\overline{S})$ by items {\it 2} and {\it 4}.

 Suppose now that $X \in I(\overline{S})$.  By definition, there is $\gamma <
\omega_2$ such that, for $i \geq \gamma,\; X \cap X_i$ is non-stationary.  Let 
$\langle T_j \mid j \in \omega_1\rangle$ be an $\omega_1$-enumeration of the
collection $\{ S_i \mid i < \gamma \}$. Then each $T_j \in I_1$ by (\Athree).
 Let $T = \nabla_{j\in\omega_1}
T_j$ be the diagonal union.  By normality of $I_1, T \in I_1$.  Hence
$X\cap T \in I_1$.  But $X\bbackslash T$ has only countable intersections with
each $T_j$ (for in fact $(X\setminus T) \cap T_j \subseteq j+1$),
 and hence, certainly, has non-stationary intersections with {\em
every} $S_i$, and is thus in $I_0$ (by (\Afour)).
As $I_0 \subseteq I_1$ (by formula (\ref{e0}) in Definition
\ref{Dideals}), $ X \in I_1$. \qed

We will prove next that if $A(\overline{S},\overline{\eta})$ holds, then
 $\overline{\eta}$ is determined, up to an $I_1$ set, as that ladder 
$\overline{\eta}$ for which $(\exists \overline{S})
 A(\overline{S},\overline{\eta})$.

\begin{lemma} 
\label{l23}
 If the first five statements hold for 
 $A(\overline{S},\overline{\eta}^1)$ and $A(\overline{T},
\overline{\eta}^2)$, then  $I(\overline{S})=  I(\overline{T})= I_1$,
and $ \{ \delta \in \omega_1\mid  
  [\eta^1_\delta] \neq^* [\eta^2_\delta] \}
\in   I_1$.
\end{lemma}
{\bf Proof}.  Define 
$S^1 =  \dom (\overline{\eta}^1\bbackslash\overline{\eta}^2)$, and 
$S^2 = \dom (\overline{\eta}^2 \bbackslash\overline{\eta}^1)$.  We claim that 
$S^1,S^2 \in I_1$.  This implies
the lemma because $S^1\cup S^2 \in I_1$ follows.
 By symmetry, it suffices to deal with only
one of these sets, for example with $S^1$.

Set $\overline{\rho} = \overline{\eta}^1 \bbackslash \overline{\eta}^2$ (so 
$S^1 = \dom(\overline{\rho}))$. Since it is
 disjoint from $\overline{\eta}^2,
\overline{\rho}$ is avoidable (by item (\At) of $A(\overline{T},
\overline{\eta}^2))$.  Yet, $\overline{\rho} \lhd \overline{\eta}^1$ , and
so, by Lemma \ref{l21} (1),
 $\dom(\overline{\rho})$ is $\overline{S}$-small, which, in view of
Lemma \ref{l21}({\it 5}), implies that $S^1 \in I_1$. \qed

Whenever $A(\overline{S},\overline{\eta})$ holds, a set of reals can be
decoded which we denote $\code(\overline{S},\overline{\eta})$.
We will encode reals (subsets of $\omega$)
 by taking subladders of $\overline{\eta}$ appropriately
chosen.  Suppose that 
$\sigma$ is a cofinal subset of order-type $\omega$ of some
$\delta < \omega_1$.  Identifying $\sigma$ with $\omega$, any $\sigma'
\subseteq \sigma$ corresponds to a subset of $\omega$.  This encoding of
reals as subsets of $\sigma$ 
is too crude, because if we take end segments of $\sigma$
and $\sigma'$ then a different real may be decoded. Since we shall be able
to recover the ladder $\overline{\eta}$ only up to finite changes we must have
a more stable decoding procedure.  So we look for a function $d$ that associates
with 
every pair $(\sigma',\sigma)$ as above some real $d(\sigma',\sigma)$ so
that:
\[ \mbox{If}\;
\sigma_1 =^* \sigma_2 \ \mbox{and}\ \sigma'_1 =^* \sigma'_2,\ \mbox{then}\
d(\sigma'_1 , \sigma_1) = d(\sigma'_2, \sigma_2).\]

The range of $d$ should be all subsets of  $\omega$, i.e., for every $\sigma$
for every $x \subseteq \omega$ there is $\sigma' \subseteq \sigma$ such that 
$d(\sigma', \sigma) = x$.  It is not difficult to find such a function $d$,
and we assume that the reader has picked one.
(For example, you may look at the intervals of $\sigma$ formed by
successive members of $\sigma'$ and take those cardinalities that appear
infinitely often.) 

Now let $\overline{\sigma}' \lhd \overline{\sigma}$ be two ladders; we say that
$(\overline{\sigma}',\overline{\sigma})$ encodes the real $r \subseteq \omega$
if, for every $\delta \in \dom (\overline{\sigma}'),\;
 d([\sigma'_\delta],[\sigma_\delta]) = r$.  We may just
write $d(\overline{\sigma}',\overline{\sigma})=r$ in such a case.

Not every pair $\overline{\sigma}'\lhd \overline{\sigma}$ encodes a real.  An
extreme case is when, for every $\delta_1 \neq \delta_2$ in
$\dom(\overline{\sigma}'),\;d(\sigma'_{\delta_1},
\sigma_{\delta_1}) \neq d(\sigma'_{\delta_2},
\sigma_{\delta_2})$.  We shall say in such a case that
 $(\overline{\sigma}',\overline{\sigma})$ are ``clearly'' not encoding.

Now we can understand the meaning of {\it A6}. If
$A(\overline{S},\overline{\eta})$ holds, we define
\[ \code(\overline{S},\overeta)=\{ r \subseteq \omega\; \mid\; r=d(\chi(S_i),
\overeta \rest S_i)\  \mbox{for some}\ i\in \omega_2\}.\]
Clearly if $r \in \code(\overline{S},\overeta)$, then $
r=d(\chi(S_i),\overeta \rest S_i)$ for an unbounded set of $i\in
\omega_2$.

\begin{lemma}
\label{l29}
If $A(\overline{S},\overeta ^1)$ and $A(\overline{T},\overeta ^2)$,
then $\code(\overline{S},\overeta ^1)= \code(\overline{T},\overeta ^2)$.
\end{lemma}
{\bf Proof}. Suppose that $r \in \code(\overline{S},\overeta ^1)$ and let
$U \subset \omega_2$ be the unbounded set of indices $i$ such that $r=
d(\chi(S_i),\overeta^1 \rest S_i)$. We must check that
for some (and hence for unboundedly many) $j\in \omega_2$, 
$r= d(\chi(T_j), \overeta ^2 \rest T_j)$
We know that $[\overeta_\delta ^1]=^* [\overeta_\delta ^2]$ 
except for an $I_1$ set, and $I_1 = I(\overline{S}) = I(\overline{T})$
(Lemma \ref{l23}). That is, if $H=\{ \delta \in \omega_2\; \mid\; 
[\overeta ^1_\delta] \not =^* [\overeta ^2_\delta] \}$, then $H \in I_1$,
 and hence $H\in I(\overline{S})$.
Thus there is an index $i \in U$ such that 
\begin{equation}
\label{Ens}
H \cap S_i\ \mbox{ is non-stationary.}
\end{equation}
That is,
\begin{enumerate}
\item $\overeta ^1 \rest S_i=^* \overeta ^2 \rest S_i$ (that is, 
$[\overeta ^1_\delta] =^* [\overeta ^2_\delta]$ for all $\delta\in
S_i$, except for a non-stationary set),
\item $d(\chi(S_i),\overeta^1\rest S_i)=r$.
\end{enumerate}
Now $\chi(S_i)$ is maximal for $S_i$ (a stationary set) and hence its domain
$S_i$ is not avoidable. So by (\Afour) of $A(\overline{T},\overeta^2)$, for
some $j\in \omega_2$, $X= S_i \cap T_j$ is stationary.
Hence $\chi(S_i)\rest X$ is maximally guessing (and non-trivial). Similarly
$\chi(T_j)\rest X$ is maximally guessing, and thus
\[ \chi (S_i) \rest X =^* \chi(T_j)\rest X\]
by the uniqueness of the maximal ladder over $X$ (namely $\chi(X)$).
Since $X \cap H$ is non-stationary (by (\ref{Ens}) above),
\[ \overeta^1 \rest X =^* \overeta ^2 \rest X\]
and thus $(\chi(T_j),\overeta^2 \rest T_j)$ encodes a real, and
$d(\chi(T_j),\overeta^2 \rest T_j)= r$.

\section{The consistency of $A(\overline{S},\overline{\eta})$}

Our aim in this section is to prove the following

\noindent
{\bf Theorem.} Assume that $2^{\aleph_0}=\aleph_1$ and
$2^{\aleph_1}=\aleph_2$. Suppose that
 $\overline{T} = \langle T_i \mid i < \omega_2\rangle$ is a
collection of $\aleph_2$ pairwise 
almost disjoint stationary subsets of $\omega_1$, and  
$\overline{\eta}$ is a ladder system such that 
\begin{description}
\item{(1)} $\overline{\eta} \rest T_i$ is guessing (but not
necessarily strongly guessing) for every $i < \omega_2$.

\item{(2)} range$(\overline{\eta}) \cap T_i$ is empty for every $i$.
\end{description}
 Then there is  a generic extension in which 
$A'(\overline{T}, \overline{\eta})$ and Martin's Axiom hold.

  The extension
is an iteration of the posets $R(\overline{\mu})$, and $P(\overline{\eta},C)$
described below. Before proving this theorem,
 however, we review some notions from
 proper forcing theory.
\subsection{Some proper forcing theory}

 This short subsection assembles some known  definitions and results on
proper forcing, such as $\alpha$-properness and $S$-properness for
a stationary set $S$.  Our notations and terms are taken (with some minor
changes) from
 Shelah's book \cite{Book} (see also \cite{Chapter}).

 Recall that if $P$ is a forcing poset and $N\prec H_\lambda$ a
 countable elementary substructure, then a condition $q\in P$ is $N$
 generic iff for every $D\in N$, dense in $P$, every extension of $q$
 is compatible with some condition in $D\cap N$. A forcing poset $P$
 is {\em proper} is for some cardinal $\lambda$, for every countable
 $N\prec H_\lambda$ such that $P\in N$, every $p\in P\cap N$ has an
 extension that is $N$ generic.

\begin{definition}[of $\alpha$-properness.]
  Let $\alpha$ be a countable ordinal.
A poset $P$ is  said to be
$\alpha$-proper iff for every  large enough cardinal 
$\lambda$, {\bf if}
 $\langle N_i \mid i \leq \alpha\rangle$ is an increasing, 
 continuous sequence of countable
elementary submodels of $H_\lambda$ such that $P \in N_0$ and 
$\langle N_j \mid j \leq i\rangle \in N_{i+1}$ for every $i < \alpha$,
{\bf then} any 
$p_0 \in P\cap N_0$ can be extended to $q \in P$ that is $N_i$-generic for
every $i\leq \alpha$.\end{definition}

\begin{definition} 
 Let $S \subseteq \omega_1$ be stationary.  A forcing poset
$P$ is $S$-proper if it is proper for structures $M$ such that $M \cap \omega_1 
\in S$. That is, $P$ is $S$-proper iff for sufficiently large $\lambda$, if $M
\prec H_\lambda$ is countable, $S,P \in M$, and $M \cap \omega_1 \in S$
 then any 
$p \in P \cap M$ can be extended to an $M$-generic condition.
\end{definition}
A stronger property is that of a poset being $S$-complete.  It means that
whenever $M \prec H_\lambda$ is countable, with $P,S \in M$, and $M \cap
\omega_1 \in S$, then every increasing and generic $\omega$-sequence of
conditions in $P \cap M$ has an upper bound in $P$.
(A sequence of conditions is {\em generic} if it intersects every
dense set of $P$ in $M$.)

The notion $(E,\alpha)$-properness is defined in Shelah 
(\cite{Book} (Chapter V). 
Just as properness is equivalent to the preservation of stationarity of
$S_{\aleph_0} (\mu)$, so is $(E,\alpha)$-properness equivalent to the
preservation of an appropriate notion of stationarity defined there.  However,
for our article, a notion of 
somewhat less generality suffices.

Let $I^\omega$ be the collection of all increasing sequences of countable
ordinals.  We write $\overline{\alpha} = \langle \alpha_i \mid
i < \omega\rangle$
 for $\overline{\alpha} \in I^\omega$.
  The club
guessing property can be regarded as a notion of non-triviality of
subsets of $I^\omega$.

\begin{definition} 
\begin{enumerate}
 \item A family $E \subseteq I^\omega$ is {\em stationary} if
for every club $C \subseteq \omega_1$ there is $\overline{\alpha} \in E$ such
that $[ \overline{\alpha} ] = \{ \alpha_i \mid i \in \omega \} \subset
C$.

\item Let $E \subseteq I^\omega$ be stationary.  We say that the poset $P$ is
$E$-proper (or $(E,\omega)$-proper, to emphasize that this notion is
related to $\omega$-properness) iff for every sufficiently large cardinal
$\lambda$, whenever $M_i \prec H_\lambda$, for $i < \omega$, are countable
with $E,P \in M_0$ and are 
such that $M_i \in M_{i+1}$ for all $i < \omega$, if

\[ \langle M_i \cap \omega_1 \mid i < \omega\rangle \in E \]
then any $p \in P \cap M_0$ can be extended to a condition which is
$M_i$-generic for every $i < \omega$.
\end{enumerate}
\end{definition}
In Shelah \cite{Book}
it is proved that the countable support iteration of posets
that are $S$-proper ($\alpha$-proper or $E$-proper) is again $S$-proper (
$\alpha$-proper or $E$-proper, respectively).
Also, if $P$ is $S$-proper ($E$-proper), then, in $V^P$, $S$ (respectively
$E$) remains stationary.

\begin{lemma} If $E \subseteq I^{\omega}$ is stationary and $P$ is an
$E$-proper poset, then $E$ remains stationary in $V^P$.
\end{lemma}
{\bf Proof.} Let $D$ be a name in $V^P$ forced by some $p\in P$ to be a
club subset of $\omega_1$. Define an $\omega_1$ sequence $\langle M_i \mid
i \in \omega_1 \rangle$ where $M_i \prec H_\lambda$ are countable with
$\langle M_i \mid i \leq j \rangle \in M_{j+1}$, and such that $p, P,D \in
M_0$. The set $C = \{ M_i \cap \omega_1 \mid i \in \omega_1 \rangle$ is
closed unbounded in $\omega_1$. Since $E$ is stationary, there is  
$\overline{\alpha} \in E$ such that $\{ \alpha_n \mid n \in \omega \}
\subset C$. Then $\alpha_n = M_{i(n)}\cap \omega_1$ and $N_n=M_{i(n)}$ is
an increasing sequence of structures with $N_n \in N_{n+1}$ and such that
$\langle N_i \cap \omega_1 \mid i < \omega \rangle \in E$. So there is an
extension $q \in P$ that is $N_i$-generic for every $i < \omega$. So for
every $i$ $q \forces \alpha _i \in D$. (Because $q$ forces that $D$ is
unbounded below $\alpha_1=\omega_1^{N_i}$.) Thus $q \forces [
\overline{\alpha}] \subseteq D$, as required.

We shall define now two subsets of $I^\omega,\; E_{\overline{\eta}}$ and
 $D_{\overline{\eta}}$, which will be used later.
\begin{definition}
\begin{enumerate}
\item  Let $\overline{\eta}$ be a
ladder system and $S = \dom (\overline{\eta})$.  Define 
$E_{\overline{\eta}}\subseteq I^\omega$ by
\[ \overline{\alpha} = \langle \alpha_i \mid i < \omega \rangle \in
 E_{\overline{\eta}} \]
iff
\begin{quotation}
\noindent
 $\overline{\alpha} \in I^\omega$ and, for 
$\delta =$ sup$\{ \alpha_i \mid
i < \omega \}$, $\delta \in S$ and 
 $\overline{\alpha}$ is an end segment of $\eta_\delta$
 (i.e., for some $k,\;
\eta_\delta(k+i) = \alpha_i$ for all $i$s).
\end{quotation}

It is obvious that $E_{\overline{\eta}}$ is stationary iff $\overline{\eta}$
is club guessing.
Thus, if 
$\ov{\eta}$ is club guessing and
$P$ is $E_{\overline{\eta}}$-proper, then $\ov{\eta}$ remains
 a guessing ladder in $V^P$.
\item
The set $D_{\overline{\mu}} \subseteq I^\omega$ ($D$ is for disjoint) 
is defined for any ladder $\overline{\mu}$
as follows:
 $\overline{\alpha} \in D_{\overline{\mu}}$ iff 
for
$\delta =$ sup$\{ \alpha_i \mid
i < \omega \}$,
either $\delta 
\not\in \dom (\overline{\mu})$ or $[\overline{\alpha}] \cap [\mu_\delta]
=^* \emptyset$.
If $\ov{\mu}$ is disjoint from $\ov{\eta}$, then $E_{\ov{\eta}} \subseteq
D_{\ov{\mu}}$. Thus, in this case, if $P$ is 
$(D_{\ov{\mu}}, \omega)$-proper, then $P$ is
$(E_{\ov{\eta}}, \omega)$-proper as well.
\end{enumerate}
\end{definition}

\subsection{The building blocks}
Two families of posets are described in this subsection: $R(\overline{\mu})$
and $P(\overline{\eta},C)$.

{\bf The poset} $R(\overline{\mu})$.  Let $\overline{\mu}$ be a
ladder over a set $S\subseteq \omega_1$. 
 The poset $R(\overline{\mu})$ introduces a generic club to
$\omega_1$ that avoids $\overline{\mu}$.  So, naturally,
$$c \in R(\overline{\mu})$$ iff
\begin{quote}
 $c \subseteq \omega_1$ is countable, closed (in
particular $\max(c) \in c)$, and for every $\delta \in S,\;[\mu_\delta]
\cap c$ is finite.
\end{quote}

 The ordering on $R(\overline{\mu})$ is end-extension.  

The cardinality
of $R(\overline{\mu})$ is the continuum.
It is clear that  $R(\overline{\mu})$ is $\omega_1 \setminus S$ complete.
A short argument is needed in order to prove that it is proper. 

Observe first that for
any condition $q\in
R(\overline{\mu})$ and dense set $D\subseteq  R(\overline{\mu})$, if
$\alpha_0=\max(q)$ then there is a closed unbounded set of ordinals $\gamma
< \omega_1$, $\alpha_0 < \gamma$, such that for every $\alpha_1$ with
$\alpha_0 <\alpha_1< \gamma$ there is an extension $q'\in D$ such that $q'
\subset \gamma$ and $\alpha_1 \in q'$ is the successor of $\alpha_0$ in
$q'$. For example, the club set can be obtained by defining a
continuous, increasing chain $\langle N_\alpha \mid \alpha \in
\omega_1\rangle$ of countable elementary substructures of some
$H_\lambda$ with $\overline{\mu}$ and the dense set $D$ in $N_0$. Then
$\langle \omega_1\cap N_\alpha \mid \alpha \in \omega_1\rangle$ is as
required.

Suppose that a
countable $M \prec H_\lambda$ and a condition 
$p_0\in R(\overline{\mu})\cap M$ are given.
We want to define an increasing, generic sequence of
conditions $p_i$ extending $p_0$ so that for $\delta = M \cap
\omega_1$, $p=\bigcup_{i\in \omega} p_i \cup
\{ \delta \}$ is a condition. 
 The case $M \cap \omega_1 \not \in S$ is
trivial and so assume that $\delta = M \cap \omega_1 \in S$. 
The problem is that we may decide infinitely
often to put $\mu_\delta(n)$ in $\bigcup_i p_i$, and then $p$ is not a
condition. 
The preliminary observation enables the construction of the sequence $p_i$ in
such a way that $p \cap [\mu_\delta] \subseteq p_0$ is finite.
The point is that when we need to extend a condition $p_i$ into a
dense set $D$, we first consider the club set formulated above (do it in the
substructure $M$) and find a limit ordinal $\gamma$ in the club that is 
in $M$.
Now
$\alpha_1 < \gamma$ is chosen so that the interval $[\alpha_1,\gamma$]
is disjoint to  $[\mu_\delta]$.
(The fact that $[\mu_\delta]$ is only  an $\omega$
sequence implies the existence of such an ordinal. 

$R(\overline{\mu})$ is not $\omega$-proper.  For suppose $M_i,\;i < \omega$
is an increasing sequence of elementary submodels such that $\alpha_i = M_i
\cap \omega_1 \in [\mu_\delta]$ for infinitely many $i$'s, where $\delta = \sup 
\{\alpha_i \mid i < \omega \}$.  Then no condition can be generic for all of 
the $M_i$s. 
However, if $\delta \not\in S$ or $\langle M_i \cap \omega_1 \mid i < \omega
\rangle$ is disjoint from $[\mu_\delta]$ (or has only a  
finite intersection) then there
is no problem in finding such a generic condition.  That is, $R(\overline{\mu})$
is $(D_{\overline{\mu}}, \omega)$-proper.
In fact, if $p_i$ is any sequence of increasing conditions where $p_i \in
M_i$ is $M_{i-1}$ generic, then $\bigcup_i p_i$ gives a condition. This
property is stronger than  $(D_{\overline{\mu}}, \omega)$-properness, but
in application we shall mix proper forcings with $R(\overline{\mu})$
forcings and hence the iteration itself is $(D_{\overline{\mu}}, \omega)$-
proper.

Hence we have the following which will be used in Lemma \ref{l33}.
\begin{lemma}
\label{l31}
Suppose that $\overline{\mu}$ is a ladder system and $A,B \subseteq \dom
(\overline{\mu})$ are such that $\overline{\mu}\rest A \cap B$ is not
guessing. Then $R(\overline{\mu}\rest A)$ is $E_{\overline{\mu} \rest
B}$-proper.
\end{lemma}
{\bf Proof}. Suppose that  $A,B
\subseteq \dom(\overline{\mu})$ are such that  $\overline{\mu}\rest A \cap
B$ is not guessing. Let $C \subseteq \omega_1$ be a club set such
that, for
every $\delta \in A \cap B$, $[\mu_\delta ] \not \subset ^* C$. Suppose that
$M_i \prec H_\lambda$ for $i<\omega$ are as in the definition of
 $E_{\overline{\mu} \rest B}$ properness and $\delta = \sup(M_i \cap
\omega_1\;\mid \; i < \omega)$. So, $E_{\overline{\mu} \rest
B}, \ R( \overline{\mu} \rest A) \in M_0$.  Hence $A, B \in M_0$ and
thus $C\in M_0$ can be assumed. Then $M_i \cap \omega_1 \in
C$ for every $i$. Since $\langle M_i \cap \omega_1\;\mid \; i <
\omega\rangle\in E_{\overline{\mu} \rest B}$, $\delta \in B$ and
$[\mu_\delta ] =^*\{ M_i \cap \omega_1\ \mid\ i < \omega \}$. Thus
$[\mu_\delta ] \subseteq ^* C$ and hence $\delta \not \in A \cap B$. So
$\delta \not \in A$ and as $R(\overline{\mu}\rest A)$ is $\omega_1
\setminus A$ complete, there is no problem in finding a condition that is
$M_i$-generic for every $i$. \qed

\noindent
{\bf The poset} $P(\overline{\eta},c)$. Let $\overline{\eta}$ be a guessing
ladder over a stationary co-stationary set $S$, such that
$$S \cap \ \range (\overline{\eta})\ \mbox{is non-stationary}.$$
(See Definition \ref{Dladder} for $ \range (\overline{\eta})$.
Then, for any club set $c \subseteq \omega_1$, the poset
$P(\overline{\eta},c)$ introduces a generic club set $D \subset \omega_1$,
such that for every $\delta \in D \cap S,\;[\eta_\delta] \subset^* c$.  This may
be viewed as forcing a club subset to the stationary set $\{ \delta \in
S \mid [\eta_\delta] \subseteq^* c\} \cup (\omega_1\bbackslash S)$.

Accordingly, we define 
$d \in P(\overline{\eta},c)$ iff $d \subseteq \omega_1$ is countable, closed
(with $\max(d) \in d$), and for every $\delta \in d \cap S,\; [\eta_\delta]
\subseteq^* c$.\\

The order is end-extension. 

 It is easy to check that any condition has
extensions to arbitrary heights (as there are no restrictions on
$\omega_1\bbackslash S$).  The cardinality of $P(\overline{\eta},c)$ is the
continuum.

$P(\overline{\eta},c)$ is not necessarily proper, because if, for $\delta = M
\cap \omega_1,\; [\eta_\delta] \not\subseteq^* c$, then no $M$-generic condition
can be found.  Still, $P(\overline{\eta},c)$ possesses two good properties
which allows its usage:

\begin{enumerate}
\item $P(\overline{\eta},c)$ is $(\omega_1\bbackslash S$)-complete (the proof
of this is obvious).

\item $P(\overline{\eta},c)$ is $(E_{\overline{\eta}}, \omega)$-proper.
($E_{\overline{\eta}}$ is stationary since $\overline{\eta}$ is
guessing.)
\end{enumerate}
We check the second property --- it is for its sake that the requirement that 
$\dom(\overline{\eta}) \cap \range (\overline{\eta}) $ is non-stationary
 was made.  So
let $\langle M_i \mid i < \omega \rangle$ be an increasing sequence of countable
elementary submodels of $H_\lambda$, with $M_i \in M_{i+1}$, and such that 
$P(\overline{\eta},c), \overline{\eta},c \in M_0$.  Denote $\delta_i = M_i
\cap \omega_1$, and $\delta = \sup \{\delta_i \mid i < \omega \}$.

The assumption is that $\langle \delta_i \mid i < \omega \rangle \in
 E_{\overline{\eta}}$, and the desired conclusion is that any $p_0 \in P \cap
M_0$ can be extended to a condition that is generic for every $M_i$.  So the
assumption is that $\delta \in S$ and $\langle \delta_i \mid i < \omega \rangle$
is an end segment of $\eta_\delta$.  Since $\dom(\overline{\eta}) \cap 
\range (\overline{\eta})$ is non-stationary, $\delta_i \not\in S$ (because$M_0$
contains a club that is disjoint from this intersection), 
and it is easy to find
(in $M_{i+1}$) an $M_i$-generic condition extending any given condition
(using the $(\omega_1 \setminus S)$-completeness).
Thus, given $p_0 \in P(\overline{\eta},c) \cap M_0$, we may construct an
increasing sequence of conditions $p_i \in M_i$, such that $p_{i+1}$ is 
$M_i$-generic.  Then $p = \{\delta\} \cup \bigcup_{i<\omega} p_i$ is in $P(\overline{\eta},c)$
because $[\eta_\delta] \subseteq^* c$ follows from the fact that 
$\delta_i \in c$ for every $i$ (as $c \in M_i$).\\

As a warm-up we shall present some simple models obtained by countable
support iteration of the posets $R(\overline{\mu})$ and
$P(\overline{\eta},c)$ just described. We assume $2^{\aleph_0}=\aleph_1$
and $2^{\aleph_1}=\aleph_2$ in the ground model.
\begin{enumerate}
\item A model in which $MA + 2^{\aleph_0}=\aleph_2 + \omega_1$ {\em is
avoidable}. This is achieved by iterating $c.c.c$ posets to obtain
Martin's Axiom, and posets of the form $R(\overline{\mu})$ (varying over
all possible ladders $\overline{\mu}$ over $\omega_1$). Countable support
is used in this iteration of proper forcing posets, and hence the final
poset is proper. 
The final poset satisfies the
$\aleph_2$-chain condition (see \cite{Book}, Chapter VIII, or
\cite{Chapter}).
The length of the iteration is $\omega_2$ so that each
possible $c.c.c$ poset of size $\aleph_1$ and each ladder $\overline{\mu}$
are taken care of at some stage. 

\item Given a guessing ladder $\overline{\eta}$ such that
$\dom(\overline{\eta}) \cap \range(\overline{\eta})=\emptyset$ a model
of $MA + 2^{\aleph_0}=\aleph_2 $ can be obtained in which
$\overline{\eta}$ is strongly guessing. 
This time posets of type $P(\overline{\eta}, c)$ 
are iterated (varying club sets $c \subseteq \omega_1$) as well as
$c.c.c.$ posets. The iteration is with countable support and of length
$\omega_2$ as before. Put $S=\dom(\overline{\eta})$.
Then $S$ is stationary (as $\overline{\eta}$ is guessing) and
co-stationary (as $S\cap \range(\overline{\eta})=\emptyset$).
Since each poset $P(\overline{\eta}, c)$ is $\omega_1 \setminus S$
complete, and each $c.c.c.$ poset is obviously $\omega_1 \setminus S$
proper, we have here an iteration of $\omega_1 \setminus S$ proper posets.
Thus the final poset itself is $\omega_1 \setminus S$ proper and $\omega_1$
is not collapsed. Moreover, since the iterands
 (both $P(\overline{\eta}, c)$ and the $c.c.c.$ posets) are
$E_{\overline{\eta}}$ proper, the final iteration is $E_{\overline{\eta}}$
proper. Hence $\overline{\eta}$ remains guessing at each stage and 
in the final extension. It
is strongly guessing since we took explicit steps to ensure this.
 
\item Now we want to combine 1 and 2. We are given a guessing ladder system
$\overline{\eta}$ defined over a stationary co-stationary set $T$, such that
$T \cap \range(\overline{\eta})=\emptyset$, and we want a generic extension
in which $\overline{\eta}$ is maximal for $\omega_1$. For the iteration,
decompose $\omega_2$ into three sets $\omega_2 = J \cup K \cup I$ of
cardinality $\aleph_2$ each. At stage $\alpha< \omega_2$ of the iteration,
supposing that $P_\alpha$ has been defined, define the poset $Q_\alpha$ in
$V^{P_\alpha}$ as follows:
\begin{enumerate}
\item If $\alpha \in J$, then $Q_\alpha$ is a $c.c.c$ poset, and the
iteration of all posets along $J$ guarantees Martin's Axiom.
\item For $\alpha\in K$, $Q_\alpha$ will be of type $R(\overline{\mu})$
where $\overline{\mu}\in V^{P_\alpha}$ is a ladder system disjoint from
$\overline{\eta}$. $R(\overline{\mu})$ is proper and it is
$(D_{\overline{\mu}}, \omega)$ proper. Hence as $E_{\overline{\eta}}
\subseteq D_{\overline{\mu}}$, $R(\overline{\mu})$ is $E_{\overline{\eta}}$
proper.

\item For $\alpha \in I$, $Q_\alpha$ will be of type 
 $P(\overline{\eta}, c)$ where $c$ is a club set in $V^{P_\alpha}$. These
posets are $\omega_2 \setminus T$ complete, and $E_{\overline{\eta}}$
proper. 
\end{enumerate}
Any of the posets along the iteration is either proper or
 $\omega_2 \setminus T$ proper
 (namely, the $P(\overline{\eta}, c)$ posets which are 
 $\omega_2 \setminus T$ complete). So the iteration itself is $\omega_1
\setminus T$ proper, and thus $\omega_1$ is not collapsed.
Moreover, the posets are $E_{\overline{\eta}}$ proper, and hence
$\overline{\eta}$ retains its guessing property in the extension.
\end{enumerate}
\subsection{The iteration scheme}
Recall that our aim is to prove the following theorem.

 \begin{theorem}  Assume $2^{\aleph_0} = \aleph_1$ and
$2^{\aleph_1} = \aleph_2$.  Suppose 

\begin{description}
\item{(1)} A sequence $\overline{T} = \langle T_i \mid i \in \omega_2\rangle$ of
pairwise almost disjoint stationary subsets of $\omega_1$.
(Almost disjoint in the sense that $T_i \cap T_j$ is non-stationary.)

\item{(2)} A ladder system $\overline{\eta}$ such that 
\begin{enumerate}
\item $\bigcup\{ T_i \mid i \in \omega_2\} \subseteq
\dom(\overline{\eta})$, and $\omega_1 \setminus \dom(\overline{\eta})$
is stationary.
\item For every $i,\;\overline{\eta} \rest T_i$ is club guessing.
 \item $\dom(\overline{\eta}) \cap$ range$(\overline{\eta}) = \emptyset$.
\end{enumerate}
\end{description}
Then there is cofinality preserving generic extension in which
 $A'(\overline{T},\overline{\eta})$ and $MA + 2^{\aleph_0} = \aleph_2$ hold.
 (The definition of $A'(\overline{T},\overline{\eta})$ is immediately
 after Definition \ref{DA}.)
\end{theorem}
{\bf Proof}. It is not difficult to get $\overline{T}$ and
$\overline{\eta}$ as in the theorem, and 
the following section contains a generic construction of
such objects.
 Here we just assume their
existence and prove the theorem.
 The generic extension is made via $P = P_{\omega_2}$, obtained
as an iteration, $\langle P_\alpha \mid \alpha \leq \omega_2\rangle$, with
countable support of posets of cardinality $\aleph_1$.  At successor stages, 
$P_{\alpha+1} \cong P_\alpha * Q_\alpha$, where $Q_\alpha \in V^{P_\alpha}$ is
one of the following three types.
\begin{description}
\item{(1)} A c.c.c. poset. (To finally obtain Martin's Axiom.)

\item{(2)} A $P(\overline{\sigma},c)$ poset, where $\overline{\sigma} \in
V^{P_\alpha}$ is a guessing ladder such that $\overline{\sigma} \lhd
 \overline{\eta}$, and $c \in V^{P_\alpha}$ is a club set.  Recall that 
$P(\overline{\sigma},c)$ introduces a generic club subset $D$ such that
$\delta \in D \cap\dom (\overline{\sigma})$ implies $\sigma_\delta
\subseteq^* c$.  We have checked that this poset is $(\omega_1\bbackslash$
 $\dom(\overline{\sigma}))$-complete, and
$(E_{\overline{\sigma}},\omega)$-proper (as $\dom(\overeta) \cap
\range(\overeta)=\emptyset$).

\item{(3)} The third type of iterated posets is $R(\overline{\mu})$ where 
$\overline{\mu} \in V^{P_\alpha}$ is a ladder system.  This forcing makes
 $\overline{\mu}$  avoidable.  We have seen that $R(\overline{\mu})$ is
proper, $(D_{\overline{\mu}}, \omega)$-proper, and $(\omega_1\bbackslash
\dom(\overline{\mu}))$-complete.
\end{description}

 Each iterated poset $Q_\alpha$ is $(\omega_1 \setminus
 \dom(\overline{\eta}))$ proper in $V^{P_\alpha}$. Hence the iteration
 itself is $(\omega_1 \setminus \dom(\overline{\eta}))$ proper, and it
 satisfies the $\alpha_2$-c.c.

 We must specify how to choose the posets $Q_\alpha$ for 
the iteration.  Every
$P_\alpha$ will have cardinality $\leq \aleph_2$ and will satisfy the
$\aleph_2$-c.c.  
When we say that a name in $V^{P_\alpha}$ satisfies property $\phi$,
we mean that it is forced by every condition in $P_\alpha$ to satisfy
$\phi$.
We say that a $V^{P_\alpha}$ name of a subset of
$\omega_1$ is {\em standard} iff it associates with every $\beta \in \omega_1$
a maximal antichain of conditions that decide whether $\beta$ is in this
subset or not.  Every subset of $\omega_1$ in $V^{P_\alpha}$ has (an
equivalent) standard name. 
For every poset $P$ of size $\aleph_2$ that satisfies the $\aleph_2$-c.c.,
 Fix an enumeration $\{ E(P,\gamma)\mid
\gamma < \omega_2\}$ of all standard names in $V^{P}$ of subsets of
$\omega_1$ and of ladder systems.  Thus any ladder or subset of $\omega_1$
 in $V^{P_\alpha}$ has a name
of the form $E(P_\alpha,\gamma)$ for some $\gamma < \omega_2$.  Fix a natural
well-ordering of the pairs $\{ \langle \alpha,\gamma\rangle \mid \alpha,
\gamma < \omega_2 \}$
that has order-type $\omega_2$. So each $\langle \alpha,\gamma\rangle$
has its ``place'' in $\omega_2$. This will serve in the choice of
$Q_\alpha$.

To define the iteration, we partition $\omega_2$ (in $V$):

\[ \omega_2 = J \cup K \cup L \cup \bigcup \{I_i \mid i < \omega_2 \} , \]
where each set in this partition has cardinality $\aleph_2$.  The type of 
$Q_\alpha$ depends on the set in this partition that contains $\alpha$.

For $\alpha \in J,\;Q_\alpha$ is a c.c.c. poset of cardinality $\aleph_1$,
 and the iteration of these  posets in $J$ shall provide Martin's Axiom.  
By now this
is so standard that no further details will be given.

For $\alpha \in K,\;Q_\alpha$ will be of type $R(\overline{\mu})$, where 
$\overline{\mu} \in V^{P_\alpha}$ is a ladder system disjoint from
 $\overline{\eta}$ (namely, $\overline{\mu}$ is a name forced by every
 condition to be a ladder-system disjoint from $\overline{\eta}$).
  The final result of iterating these posets along $K$ is
that, every $\overline{\mu} \in V^P$ disjoint from $\overline{\eta}$ 
is avoidable in $V^P$.  
Thus property (\At) of $A'(\overline{T}, \overline{\eta})$ can be assured.

Before going on, let's discuss the problem involved in the direct approach
to obtain (\Afive) and why we do not get $A(\overline{T},\overeta)$ but rather
$A'(\overline{T},\overeta)$ (namely $A(\overline{S},\overeta)$
where $S_i \subseteq T_i$). A possible approach to (\Afive) is to consider each
possible ladder $\overro \lhd \overeta$ such that $X=\dom(\overro)$ is not
$\overline{T}$-small, and to find for this $\overro$ some $i \in \omega_2$
such that $X \cap T_i$ is stationary. Then, if possible, to transform
$\overro \rest X \cap T_i$ into a maximally guessing ladder. For this to
have any chance, it must be the case that $\overro \rest T_i$ is guessing.
Yet it is possible that $\overro \rest T_i$ is non-guessing for every $i$.
In this case we must shrink the $T_i$'s so as  to make $X$ 
$\overline{S}$-small. This shows the need for defining subsets $S_i
\subseteq T_i$. But now (\Afour) causes a problem because, if $X\subseteq
\omega_1$ is such that in $V^P$ $X \cap S_i$ is non-stationary for every
$i<\omega_2$, then we must be able to identify this $X$ at some
intermediary stage of the iteration so as to make $\overeta \rest X$
avoidable. Yet, as the $S_i$ are not yet all defined in any intermediate
stage, it is not clear how to identify these $X$'s.

We describe now in general terms
 how the sets $I_i$ from the partition will be used in the
iteration.
For every $ i<\omega_2$ 
let $\alpha(i)$ be the first ordinal in $I_i$.  A stationary
subset $S_i \subseteq T_i$ and a guessing 
ladder $\overline{\sigma}^i$ over $S_i$ 
will be
defined in $V^{P_{\alpha(i)}}$.  The iteration  of the posets
$Q_\alpha$ for $\alpha \in I_i$ will  make $\overline{\sigma}^i$ maximal 
for $T_i$, and will achieve (\Athree)
by establishing $\chi(T_i) = \overline{\sigma}^i$.
Finally, in $V^P$, $A(\ov{S}, \ov{\eta})$
 will hold for $\ov{S} = \langle S_i \mid i \in \omega_2 \rangle$.

To define $S_i$ and $\overline{\sigma}^i$ we assume a function,
$\rho$,  which assigns to any $\alpha$ of the form $\alpha =
\alpha(i)$ a name $\rho(\alpha)\in V^{P_\alpha}$ that is one
of the following.
\begin{enumerate}
\item  If $i$ is an even ordinal then $\rho(\alpha(i))$ is 
a name of a real in $V^{P_\alpha}$. The complete definition of $\rho$
is given in the following section where it is used to define the
encoding of the well-ordering of reals. Here we only assume that
$\rho(\alpha)$ is defined.

\item If $i$ is an odd ordinal, then $\rho(\alpha(i))$ is  determined
as  a name of the ladder
system $\rho(\alpha)$, defined as follows in $V^{P_\alpha}$. 
With respect to the well-ordering of names in  $V^{P_\alpha}$,
$\rho(\alpha)$ is the least
ladder $\overline{\rho} \lhd \overline{\eta}$ that is not of the form
$\rho(\alpha')$ for $\alpha' < \alpha$, and is such that for $D =
 \dom(\overline{\rho})$
\[ \overline{\eta} \rest D \cap T_i \;\mbox{is guessing}. \]
\end{enumerate}

Suppose that  $\alpha = \alpha(i)$. Instead of defining the names
$S_i$ and $Q_\alpha$ directly in $V^{P_\alpha}$, we let 
$G \subseteq P_\alpha$ be $V$-generic and we shall describe the
interpretations of $S_i$ and $Q_\alpha$. 
We will later 
see (Lemma \ref{l33})
 that $\overline{\eta} \rest T_i$ remains guessing in
$V[G]$. In $V[G]$, collect all sets $X \subseteq \omega_1$ 
such that 
\begin{enumerate}
\item a standard name of $X$ appeared
before $\alpha$ in the well-ordering of the names (i.e., for some 
$\alpha' \leq \alpha$ and $\gamma< \omega_2$,
 $\langle \alpha',\gamma \rangle$ is
placed before $\alpha$ in the well-ordering of $\omega_2 \times \omega_2$,
and $E(P_{\alpha'}, \gamma)$, 
the $\gamma$th name in $V^{P_{\alpha'}}$, gives $X$),  and 
\item $X$ is such that 
$\overline{\eta} \rest (T_i \cap X)$ is not guessing (i.e., 
$T_i \cap X \in I_{\overline{\eta}})$.
\end{enumerate} 
Let $\langle X_\xi \mid \xi < \omega_1\rangle$ 
be an enumeration of these sets.   Take their diagonal
union
\begin{equation}
\label{e2}
A = \nabla_{\xi\in\omega_1} (X_\xi \cap T_i) .
\end{equation}
 Then
 $A \in I_{\overline{\eta}}$. 
Since $\oveta \rest T_i$ is guessing in $V[G]$, 
$T'_i = T_i\bbackslash A \not\in I_{\overline{\eta}}$  (that is, 
$\overline{\eta} \rest T'_i$ is guessing). 
(The reason for this specific definition of $A$ and $T'_i$ will only be
apparent in the proof of item (\Afive) in $V^P$.)

Now $\rho(\alpha)$ is either a real or a ladder system in $V[G]$.
Accordingly the definition of $S_i$ and $\oversigma ^i$ is split in two.
Suppose that $\rho(\alpha)$ is a real $r\subseteq \omega$ in $V[G]$. We
want to encode $r$. Define $S_i=T'_i$, and let $\oversigma ^i \lhd \overeta
\rest S_i$ be a ladder system over $S_i$ such that
\[ d(\oversigma ^i, \overeta \rest S_i)= r.\]
Since $\overeta \rest S_i$ is guessing and 
$\oversigma ^i \lhd \overeta \rest S_i$
has domain $S_i$, $\oversigma ^i$ is also guessing.

Suppose next that $\alpha = \alpha(i)$ for $i$ an odd ordinal and 
$\overro=\rho(\alpha)$ is (in $V[G]$) a ladder over
$X=\dom(\overro)$ (such that $\overro \lhd \overeta$, and 
$\overeta \rest X \cap T_i$ is guessing).
 Then $\overline{\eta} \rest X \cap T'_i$ is guessing, because
$T_i \setminus T'_i \in I_{\ov{\eta}}$. 
It follows that 
$\overline{\rho} \rest X \cap T'_i$ is guessing as well, because
$\overline{\rho} \lhd \overline{\eta}$ and $X \cap T'_i \subseteq
 \dom(\overline{\rho})$.
 In this case define 

$$S_i = X \cap T'_i,$$
and define
$$ \overline{\sigma}^i \lhd  \overline{\rho} \rest S_i\
\mbox{so that}\ (\oversigma ^i, \overeta \rest S_i)\ \mbox{is clearly not
encoding}.$$

The iteration along $I_i$ builds up the properties of $\overline{\sigma}^i$
and establishes $\chi(T_i) = \overline{\sigma}^i$ in $V^P$.  For this, the
posets $Q_\xi$, for $\xi \in I_i$, are of two types:
\begin{description}
\item{(1)} $P(\overline{\sigma}^i,c)$, where $c$ ``runs'' over all possible
clubs.  This ensures that $\overline{\sigma}^i$ becomes strongly club guessing
in $V^P$.  To enable the use of $P(\overline{\sigma}^i,c)$ we rely on the
assumption, proved later to hold, that $\overline{\sigma}^i$ remains club
guessing at each stage.

\item{(2)} $R(\overline{\mu})$, where $\overline{\mu}$ ``runs'' over all
possible ladders over $S_i$ that are disjoint from $\overline{\sigma}^i$.  This
ensures the maximality of $\overline{\sigma}^i$.
\end{description}
To satisfy item (\Afour) (in the definition of $A(\ov{S},\ov{\eta})$),
 every $X$ must be made avoidable whenever all the
intersections $X \cap S_i$ are non-stationary.  It suffices to show in such a
case that the ladder $\overline{\eta} \rest X$ is avoidable to
conclude that $X$ is avoidable, because any ladder disjoint from
 $\overline{\eta}$ is necessarily avoidable in $V^P$.  It is the iteration
along $L$ that achieves this, by forcing with posets of type $R(\overline{\eta}
\rest X)$ as follows. 

 Given $\zeta \in L$ and a generic filter $G \subseteq
P_\zeta$, we will define $Q_\zeta$ in $V[G]$.  For $i< \omega_2$
such that  $\alpha(i) < \zeta$,
the sets $S_i$ have been defined.  For {\em every} $i < \omega_2$ define

\[ S^*_i = \left\{ \begin{array}{ll}
S_i& \;\mbox{if} \;\alpha(i) < \zeta \\
T_i& \;\mbox{ otherwise}\end{array} \right. \]
Using the well-ordering of standard names,
take the least set $X \subseteq \omega_1$ (if there is one) that was not taken
before at a stage in $L$, such that

\begin{equation}
\label{st1}
 \forall i < \omega_2\;\;\overline{\eta} \rest X
\cap S^*_i \;\;\mbox{is not guessing}. 
\end{equation}
Then define $Q_\zeta$ to be $R(\overline{\eta} \rest X)$ (or a
trivial poset if no such $X$ exists).

This ends the definition of the iteration, but
it is not yet clear why items (\Afour) and (\Afive)  hold in $V^P$.
To prove (\Afour) we shall first prove that
if $$(\forall i < \omega_2) X \cap S_i \ \mbox{is non-stationary in }
V^P,$$ then (\ref{st1})
 holds at some stage
$V^{P_\xi}$, $\xi \in L$, and hence $\overline{\eta} \rest X$ is avoidable
in the next step of the iteration.
To see that this is
indeed the case, we need the following pivotal observation.

\begin{lemma} 
\label{lemma38}
 Suppose $G \subseteq P_{\omega_2}$ is $V$-generic.  If $\gamma <
\omega_2$ and $\overline{\rho} \in V[G \rest \gamma]$ is a 
ladder over
$X$ such that, in $V[G \rest \gamma ]$
$\overro \lhd \overeta$ and 

\begin{equation}
\label{st2}
 |\{i \in \omega_2 \mid \overline{\eta} \rest X \cap
T_i\;\mbox{is guessing}\}| = \aleph_2 . 
\end{equation}
Then there is $i$ such that $S_i \subseteq X \cap T_i$ and
$\overline{\sigma}^i \lhd \overline{\rho}$.\end{lemma}
{\bf Proof}.  The proof of this lemma depends on the fact that for any $i$
(with $\alpha(i) \geq \gamma$) such that $\overline{\eta} \rest X
\cap T_i$ is guessing in $V[G \rest \gamma],\;\overline{\eta}
\rest X \cap T_i$ remains guessing in $V[G \rest \alpha(i)]$ as well.  Thus, as the turn of $\overline{\rho}$ cannot be delayed
$\omega_2$ many times, at some stage $\alpha = \alpha(i),\;\overline{\rho} =
\rho(\alpha)$ holds, and then $S_i = X \cap T'_i$ and
$\overline{\sigma}^i\lhd
\ov{\rho}
\rest S_i$ were defined in $V[G \rest \alpha].$ 
\qed

Now we can prove item (\Afour) in $V[G]$.  For this, 
let $X \subseteq \omega_1$ be such that
$(\forall i < \omega_2)\;X \cap S_i$ is non-stationary.  We will show that, at
some stage $\zeta \in L,$ the poset $R(\overline{\eta} \rest X)$ was 
taken as
$Q_\zeta$. 
 If, for some $\gamma < \omega_2$,  (\ref{st2}) of Lemma  \ref{lemma38}
holds in $V[G \rest
\gamma]$, then $S_i \subseteq
X$ contradicts the fact that $S_i$ is stationary.  Hence  formula
 (\ref{st2}) never holds,
and for $\gamma$ such that $X \in V[G \rest \gamma]$ there are only
boundedly many $j$s for which $\overline{\eta} \rest X \cap T_j$ is
guessing. So let $\gamma < j_0 < \omega_2$ be such that if $\ov{\eta} \rest X \cap
T_j$ is guessing, then $j < j_0$.
  Since, in $V[G]$, $X \cap S_j$ is non-stationary for every $j < \omega_2$,
there is a stage $j_1$ such that for every $j < j_0$,
$X \cap S_j$ is non-stationary in $V[G \rest j_1]$. Thus, for
$\zeta \geq j_1$, 
 in $V[G \rest \zeta]$, 
for every $i <
\omega_2$,
if $i< j_0$ then $S_i$ is defined (that is, $\alpha(i) \leq
\zeta)$  and $ X \cap S_i$ is non-stationary, 
and hence $\overline{\eta} \rest X \cap S_i$ is  non-guessing, and
if $i\geq j_0$, then $\overline{\eta} \rest X \cap T_i$ is
non-guessing.  But this is exactly the condition required at stages 
$\zeta \in L$ to force with $R(\overline{\eta} \rest X )$.

Finally, we turn to prove item (\Afive).  So let $ \overline{\rho} \lhd
 \overline{\eta}$ with $X = \dom(\ovrho)$ 
 be given in the generic extension $V[G]$.  Then for some
$\gamma < \omega_2$,
$\overline{\rho} \in V[G\rest \gamma]$, and a name of $X$
appeared before $\gamma$ in the well-ordering of names.

{\bf Case 1:}  In $V[G\rest \gamma]$:  There is $i_0 < \omega_2$
such that, for every $i \geq i_0$, $\alpha(i) > \gamma$ and

\[ \overline{\eta} \rest X \cap T_i\;\mbox{is not guessing}. \]
Then, in defining $S_i$ for $i \geq i_0$, $\alpha(i) > \gamma$,  and the
set $X$ appears as some $X_\xi$ (in equation (\ref{e2})), 
and hence $X \cap S_i$ is at most
countable (it is included in $\xi+1$).  Thus, in Case 1, $X \cap S_i$ is
non-stationary (and even countable) for a co-bounded set of indices.  That is,
$X$ is $\overline{S}$-small. (It is for this argument that, in defining
$S_i$, we asked $S_i \cap A = \emptyset$)

{\bf Case 2}:  Not Case 1.  Hence (\ref{st2}) holds in $V[G\rest
\gamma]$.  So, by Lemma \ref{lemma38} there is $i$ such that $S_i \subseteq X
\cap T_i$ and $\overline{\sigma}^i \lhd \overline{\rho}$, which establishes 
(\Afive).  \qed

Our proof relied on preservation claims that some ladders retain their
guessing property, and we intend now to prove these claims.  First, 
set $T = \cup T_i$. Then $\omega_1 \setminus T$ is stationary by
assumption, and 
all the posets used are $(\omega_1\bbackslash
T)$-proper.  (The c.c.c. posets are certainly proper.  The
$P(\overline{\sigma},c)$ posets (defined for $\overline{\sigma}\lhd
 \overline{\eta})$ are $(\omega_1\bbackslash
\dom (\overline{\sigma}))$-complete, and hence $(\omega_1\bbackslash T)$-proper.
 The $R(\overline{\mu})$ posets are proper.)  This secures the preservation of
$\aleph_1$.

\begin{lemma}
\label{l33}
$\overline{\eta} \rest T_i$ remains guessing
in $V^{P_\alpha}$ for $\alpha = \alpha(i)$.  
\end{lemma}
{\bf Proof.} This follows from the fact that
the posets iterated at stages $\zeta < \alpha$ are all $(E_{ \overline{\eta}
\rest T_i}, \omega)$-proper: 
\begin{enumerate}
\item The c.c.c. posets are always
$\omega$-proper.  
\item The $R(\ov{\mu})$ posets iterated at stages in $K$ are defined for
$\ov{\mu}$'s that are disjoint from $\oveta$. In such a case $E_{\oveta}
\subseteq D_{\ov{\mu}}$. But we remarked that $R(\ov{\mu})$ is
$(D_{\ov{\mu}}, \omega)$-proper.
\item The $P(\overline{\sigma},c)$ posets introduced for
$\zeta < \alpha$ are defined along $I_j$
only for $j$s such that $\alpha(j) < \alpha$, and thence for
$\ov{\sigma}$'s such that $\overline{\sigma}
\lhd \overline{\eta} \rest T_j$, implying the
$(E_{\overline{\eta}\rest T_i},\omega)$-properness. 
(Since 
 $T_j \cap T_i$ is non-stationary,
and $\range(\oveta) \cap T_j = \emptyset$.)
\item  The $R(\overline{\mu})$ posets defined along $I_j$ 
for $\alpha(j) < \alpha$ are
defined for ladders $\overline{\mu}$ over $T_j$. As $T_j$ is 
almost disjoint from
$T_i$,  these $R(\ov{\mu})$'s are $(E_{\overline{\eta}\rest T_i},\omega)$-proper.
\item 
The $R(\overline{\eta}\rest X) $ posets defined for $\zeta \in
L,\ \zeta < \alpha$, are such
that $\overline{\eta} \rest X \cap T_i$ is non-guessing and the
poset is thence $(E_{\overline{\eta}\rest T_i},\omega)$-proper
(by Lemma \ref{l31}).
\end{enumerate}

Then, we must  also 
show that the guessing ladder $\overline{\sigma}^i \lhd 
\overline{\eta}\rest S_i$ defined in $V^{P_{\alpha(i)}}$ remains
guessing at every stage in $I_i$ (and thus the posets
 $P(\overline{\sigma}^i,c)$ can be applied).  This is basically the same
proof, done in $V^{P_{\alpha(i)}}$ for the quotient poset $P/P_{\alpha(i)}$
which is again a countable support iteration of posets as above that are
$E_{\overline{\sigma}^i}$ proper.

\section{The $\Sigma^{2[\aleph_1]}$ well-ordering}

 The main theorem, Theorem \ref{main}, is proved in this section.  So 
$2^{\aleph_0} = \aleph_1$ and $2^{\aleph_1} = \aleph_2$ are assumed in 
the ground model $V$.  We need
a sequence $\overline{T}$ of pairwise almost disjoint stationary sets, and a
guessing ladder system $\overline{\eta}$; and we are going to define them
first.  

Since we want to describe the $\aleph_2$ stationary sets in the language 
$\Sigma^{2[\aleph_1]}$, we need a compact form of generation for such sets. 
This is provided by the following definition.

\begin{definition} Let 
 $^{\leq \alpha}\{ 0,1 \}$ denote the set of all functions $f : \beta
\to \{ 0 , 1 \}$ for $\beta \leq \alpha$. Ordered by function
extension, this forms a tree. Define $^{< \alpha}\{ 0,1
\}$ similarly. 

A {\em stationarity tree} is a  subtree $T\subseteq ^{<\omega_1} \{
0, 1 \}$ of cardinality 
$\aleph_1$ 
 such that:
\begin{enumerate}
\item If $f,g \in T$ and $\alpha$ are such that $f \rest \alpha = 
g \rest \alpha$ but $f(\alpha) \neq g(\alpha)$, then $f^{-1} \{1\}
\cap g^{-1} \{1\} \subseteq \alpha$.
\item  $T$ has $\aleph_2$ branches of length $\omega_1$ 
and each gives a stationary set
(that is, the union of the nodes along any
$\omega_1$-branch 
forms a function $f$ and $f^{-1} \{1\}$ is a stationary subset 
of $\omega_1$). It follows from item {\it 1} that the intersection of
any pair of these stationary sets is countable.
Thus the branches of $T$ give $\aleph_2$ pairwise disjoint 
stationary sets enumerated as
$T_i$ for $i < \omega_2$.  
\item 
We also require that $U =
\cup_{i<\omega_2}\; T_i$ is a co-stationary set.  
\end{enumerate}
\end{definition}
  The poset $S$, defined below, will  produce a stationarity tree by forcing.

Conditions in $S$ 
will be countable trees,
together with countable information on the family of branches.
 Define $p \in S$ iff
 $p = (T_p,f_p)$ where:
\begin{enumerate}
\item For some countable ordinal $\beta$ (called the ``height'' of
$p$) $T_p\subseteq ^{\leq \beta}\{ 0, 1 \}$
is a countable tree of functions
 ordered by inclusion, and satisfying
property {\it 1}, and such that $T_p \cap\; ^{\beta}\{ 0, 1 \} \not =
\emptyset$.
\item  $f_p$ is a countable (partial) map defined on $\omega_2$
that assigns to $\zeta$ in its domain a node
 $f_p(\zeta)\in  T_p \cap \; ^\beta \{ 0 , 1 \}$.
($\dom(f_p)$ is called the
``domain'' of $p$.)
\end{enumerate}

The extension relation 
$p_2 \geq p_1$ on $S$ is defined by requiring that $T_{p_1} = T_{p_2} 
\cap\; ^{\leq \alpha_1} \{ 0,1 \}$
 where $\alpha_1 = $ height $(T_{p_1})$, and that 
$f_{p_2}(\zeta)$ extends $f_{p_1}(\zeta)$ for every $\zeta \in \dom (f_{p_1})$.

If $p_1,p_2 \in S$ are such that $T_{p_1} = T_{p_2}$, and $f_{p_1}$ agrees
with $f_{p_2}$ on the intersection of the domains, then $p_1$ and $p_2$ are
compatible.  Hence, CH implies the $\aleph_2$-c.c. for $S$.

It is not difficult to prove that every condition has arbitrarily high
extensions, and that for every $\xi\in \omega_2$ the set of conditions $p$
with $\xi \in \dom(f_p)$ is dense in $S$. Clearly, 
$S$ is countably closed. If $\langle p_n \mid n \in \omega \rangle$ is
an increasing sequence of conditions, let $p = \sup\{ p_n \mid n <
\omega \}$ be defined as follows. $T_p= \bigcup \{ T_{p_n} \mid n \in
\omega\} \cup \{ f_p(\xi) \mid \xi \in \dom(f_{p_n})\ \mbox{ for some
} n \in \omega \}$ where $f_p(\xi) = \bigcup \{ f_{p_n}(\xi) \mid n\in
\omega\; \&\; \xi \in \dom (p_n)\}$. That is, if $\alpha$ is the height of
$p$, then $T_p \cap\; ^\alpha \{0,1\}$ consists only of the functions
$f_p(\xi)$ for $\xi \in \dom(p)$. If $G \subseteq S$ is
$V$-generic, define $T = T(G) = \bigcup \{ T_p \mid p \in G\}$.  Then $T$ is a
stationarity tree.  
  Define 
$f(\zeta) = \bigcup \{f_p(\zeta) \mid p \in G \}$.  Then $f(\zeta)$ is
an $\omega_1$-branch of
$T$, and for $\zeta_1 \neq \zeta_2\;f(\zeta_1) \neq f(\zeta_2)$.  Thus $T$ has
$\aleph_2$ many $\omega_1$-branches.
The fact that every $f(\xi)$
gives a stationary set requires a simple density argument.
We will check now the following:\\
\begin{claim}
\label{C42} Any $\omega_1$-branch of $T$ in $V[G]$ is some $f(\zeta)$.  
\end{claim}
{\bf Proof:} Suppose, toward a contradiction, that $p$ forces
that  $\tau$
is a branch of $T $ which is not  $f(\zeta)$ for any 
$\zeta \in \omega_2$. Observe first that since $S$ is $\sigma$-closed, any
condition $p$ can be extended to a condition $q$ that describes $\tau$ up to
height$(p)$. Then, every condition $p$ and $\xi \in \dom(p)$
have an extension $q$ such that  $f_q(\zeta)$ diverges from the
value of $\tau$ determined by $q$. Repeating this procedure $\omega^2$
times, we finally get an extension $q$ of $p$ with $\delta$ = height$(q)$
limit, and such that $q$ determines $\tau$ as a branch of $T_q$ of height
$\delta$ which is different from each of the branches $f_q(\zeta)$. 
Since $T_q \cap\; ^\delta\{0,1\}$ consists only of the branches of the
form $f_q(\xi)$, the branch of $\tau$ is not in $T_q$. 
Then $q$ forces $\tau \subseteq T \rest \delta$.
\qed

We denote with 
$\overline{T} = \langle T_i \mid i < \omega_2\rangle$ the collection of
stationary sets thus obtained from the branches $f(\zeta)$ of $T$.  Let $U =
\cup_{i<\omega_2} T_i$.  

It is not difficult to show that $\omega_1\setminus U$ is also
stationary. If $C$ is a name of a closed unbounded subset of
$\omega_1$, find a countable $M\prec H_\lambda$ with $C\in M$, and
define an $M$-generic condition $p$ that puts $\delta=M\cap \omega_1$
in $\omega_1\setminus U$.

Next, we obtain a ladder system $\overline{\eta}$ over $U$ such that
 range$(\overline{\eta}) \cap U = \emptyset$, and $\overline{\eta}
 \rest T_i$ is guessing for every $i$.  It is possible to get this 
$\overline{\eta}$ by forcing with the natural (countable) conditions.  This
forcing notion is countably closed, and, assuming CH, it has cardinality 
$\aleph_1$.

Now comes the main stage of the iteration.

Using the construction of the previous section we 
obtain an extension in which $A'(\overline{T},
\overline{\eta}) + MA + 2^{\aleph_0} = \aleph_2$ hold, and such that for every 
$i < \omega_2$ either  $(\overline{\sigma}_i, \overline{\eta})$ is clearly not encoding
(where $\overline{\sigma}_i=\chi(T_i)$ is the maximal ladder for $T_i$), 
or else it encodes a real $r_i$, and in that case $r_i = r_j$ for $\aleph_2$
many $j$s (any encoded real is encoded unboundedly often).  The set of encoded
reals, ${\cal E} = \{r_i  \mid i < \omega_2\}$, is (in some natural encoding of pairs) our
well-ordering of the reals.  
We must prove that this well-ordering is $\Sigma^{2[\aleph_1]}$.
After the extension, Claim \ref{C42} may no longer be true because new
branches were added to $T$. However, the stationary sets $T_i$ are
$\Sigma^{2[\aleph_1]}$ definable. They are exactly those stationary sets
$X$ obtained from a branch of $T$ and such that $X$ is not avoidable (any
$\omega_1$-branch 
of $T$ that is not one of the original $f(\xi)$ branches is almost
disjoint to any original branch and hence by (\Afour) its stationary set is
avoidable).

We describe the $\Sigma^{2[\aleph_1]}$ formula $\psi(x)$ that decodes
this well-ordering $(\psi(x)$ iff $x \in \CE)$.  First consider the
$\Sigma^{2[\aleph_1]}$ formula $\varphi(T_0, \overline{\eta}_0)$(with
class variables $T_0$ and $\ov{\eta}_0$) which says that
$T_0$ is a stationarity 
tree, and $\overline{\eta}_0$ is a ladder system such that 
$A'(\overline{T}_0,\overline{\eta}_0)$ holds (where $\overline{T}_0$ is the
collection of non-avoidable 
stationary sets derived from the branches of $T_0$).  
  (The statement ``there are $\aleph_2$ indices such
that\ldots'' can be expressed by saying ``there is no $\aleph_1$-class
containing all the indices such that\ldots'').

  This enables a $\Sigma^{2[\aleph_1]}$ rendering
of the formula $x \in \CE$:

\[ \begin{array} {lll}
\psi(x) &\equiv &\;\mbox{there exists}\; T^0\;\mbox{and}\;\overline{\eta}^0
\;
\mbox{such that}\;\varphi(T^0,\overline{\eta}^0)\; \mbox{and}\\ &&
x \in \code(S^0,\overline{\eta}^0)
\end{array} \]

Clearly, $\psi(x)$ holds for every $x \in \CE$ (by virtue of the ``real''
$T$ and $\ov{\eta}$), and we must also prove
that if $\psi(x)$ then $x \in \CE$.
But this follow from Lemma \ref{l29}.


\begin{thebibliography}{30}

\bibitem{Chapter} U. Abraham, Proper Forcing, in Foreman, Kanamori,
and Magidor, eds, {\it Handbook of Set Theory}.
\bibitem{AS1} U. Abraham and S. Shelah,  A $\Delta_2^2$ well-order of the
reals and incompactness of $L(Q^{\rm MM})$.
 Annals of Pure and Applied Logic 
{\bf 59} (1993) 1-32.

\bibitem{AS2} U. Abraham and S. Shelah,  Martin's Axiom and $\Delta^2_1$
well-ordering of the reals. 
Archive for Mathematical Logic (1996) 35, 287--298.


\bibitem{Harrington} L. Harrington,  Long projective well-orderings. Annals of
Mathematical Logic, 12 (1977), 1--24.

\bibitem{JensenSol} R. B. Jensen and R. M. Solovay,  {Some applications of
almost disjoint sets}, in:  Y. Bar-Hillel, ed., 
{\it Mathematical Logic and Foundation of
Set Theory} (North-Holland, Amsterdam, 1970) 84-104.

\bibitem{Book} S. Shelah, {\em Proper and improper forcing}, 2nd edition,
Springer 1998. 

\bibitem{ShelahWoodin} S. Shelah and H. Woodin, { Large cardinal imply every
reasonably definable set is measurable}, Israel J. Math. {\bf 70} (1990) 
381-394.

\bibitem{Solovay} R. M. Solovay, a paper to be published in the Archive.

\bibitem{Woodin} H. Woodin, {\it Large Cardinals and Determinacy}, in
preparation.

\end{thebibliography}
\end{document}